\documentclass[a4paper,reqno,11pt, oneside]{amsart}
\usepackage{amssymb}
\usepackage{amstext}
\usepackage{amsmath}
\usepackage{amscd}
\usepackage{amsthm}
\usepackage{amsfonts}
\usepackage{enumerate}
\usepackage{graphicx}
\usepackage{bm}
\usepackage{ascmac} 
\usepackage{fancybox}
\usepackage{cancel}
\usepackage{cases}
\usepackage{ulem}
\usepackage{caption}

\newtheorem{df}{Definition}[section]
\newtheorem{thm}[df]{Theorem}

\theoremstyle{remark}
\newtheorem{rem}[df]{Remark}

\newcommand{\ZZ}{\mathbb{Z}}

\newcommand{\CC}{\mathbb{C}}
\newcommand{\RR}{\mathbb{R}}

\newcommand{\Ac}{\mathcal{A}}
\newcommand{\Bc}{\mathcal{B}}
\newcommand{\Cc}{\mathcal{C}}

\newtheorem{ex}[df]{Example}

\newcommand{\vanish}[1]{}

\begin{document}

\title[Roots of Ehrhart polynomials for $C_d^*$]{The distribution of roots of Ehrhart polynomials for the dual of root polytopes of type $C$}
\author{Akihiro Higashitani}
\author{Yumi Yamada}

\address[A. Higashitani]{Department of Pure and Applied Mathematics, Graduate School of Information Science and Technology, Osaka University, Suita, Osaka 565-0871, Japan}
\email{higashitani@ist.osaka-u.ac.jp}
\address[Y. Yamada]{Faculty of Engineering, Information and Systems, University of Tsukuba, Tsukuba, Ibaraki 305-8573, Japan}
\email{s1830124@s.tsukuba.ac.jp}

\thanks{Data availability: Not applicable.}

\maketitle

\begin{abstract}
In this paper, we study the Ehrhart polynomial of the dual of the root polytope of type $C$ of dimension $d$, denoted by $C_d^*$. 
We prove that the roots of the Ehrhart polynomial of $C_d^*$ have the same real part $-1/2$, 
and we also prove that the Ehrhart polynomials of $C_d^*$ for $d=1,2,\cdots$ have the interlacing property. 
\end{abstract}

\medskip

\section{Introduction}
A lattice polytope is a convex polytope in $\RR^d$ all of whose vertices are in the lattice $\ZZ^d$. 
For a $d$-dimensional lattice polytope $P\subset\RR^d$ there exists a polynomial in $k$ of degree $d$, 
called the \textit{Ehrhart polynomial} of $P$, that counts the total lattice points in $kP$ for all positive integers $k$, where $kP=\{ k \alpha : \alpha \in P\}$. 
Note that we write $E_P(k)$ or $|kP\cap\ZZ^d|$ for the Ehrhart polynomial of $P$. 
Ehrhart polynomials of lattice polytopes appear in many branches of mathematics and have been investigated from several points of view, 
e.g., enumerative combinatorics, toric geometry, commutative algebra, and so on. For the introduction to Ehrhart polynomials, see \cite{BR}. 

A lattice polytope $P \subset \RR^d$ whose dual polytope 
\[P^{*}:=\{u\in\RR^{d}:\langle u,v\rangle\le1\text{ for all } v\in P\}\]
is also a lattice polytope is called \textit{reflexive}, where $\langle \cdot, \cdot \rangle$ denotes the usual inner product. 
If $P$ is reflexive, then $P^*$ is also reflexive, and $P^{\circ}\cap\ZZ^d=\{0\}$, where $P^{\circ}$
%=P\backslash\partial P$
denotes the strict interior of $P$.
It is known what a lattice polytope $P$ of dimension $d$ is reflexive is equivalent to 
what the Ehrhart polynomial of $P$ satisfies the functional equation: $E_P(k)=(-1)^dE_P(-k-1)$ (see, e.g., \cite[Section 3]{BHW}). 
This implies that the roots of the Ehrhart polynomial of $P$ distribute symmetrically with respect to the line of $\mathrm{Re}(z)=-1/2$. 
%The distributions of the roots of the Ehrhart polynomials of reflexive polytopes are studied by many researchers (cf. \cite{BHW,HHK}).
The vertical line $\mathrm{Re}(z)=-1/2$ is called the \textit{canonical line}. 
It is natural to ask what kinds of reflexive polytopes satisfy that all roots of their Ehrhart polynomials lie on the canonical line. 
Such reflexive polytopes
%(canonical line hypothesis) 
are called a \textit{CL-polytope} (\cite[Definition 3]{HHK}) and CL-polytopes have been intensively studied. See, e.g., \cite{G, HHK, MHNOH, V}. 

It is hard to show that a reflexive polytope is a CL-polytope, but the theory of interlacing polynomials, which is explained below, plays an essential role for the proof. 
For two polynomials $f$ and $g$ of degree $d$ and $d-1$, respectively, 
assume that all the roots of $f$ and $g$ lie on the same line $L=\beta+\RR\gamma$ on the complex plane $\CC$ with $\beta,\gamma\in\CC$. 
Let $\beta+t_{1}\gamma,\cdots,\beta+t_{d}\gamma$ and $\beta+s_{1}\gamma,\cdots,\beta+s_{d-1}\gamma$ be the roots of $f$ and $g$, respectively.  
We say $f$ is \textit{$L$-interlaced} by $g$ if we have 
\[t_{1}\le s_{1}\le t_{2}\le\cdots\le t_{d-1}\le s_{d-1}\le t_{d}.\]
For a series of polynomials $\{f_{i}\}_i$ of degree $i$, where all the roots of each $f_i$ lie on the same line $L$, 
we say the series of polynomials satisfies the \textit{interlacing property} if $f_{i+1}$ is $L$-interlaced by $f_{i}$ for all $i$.
This interlacing property is the key in \cite{HKM} for the proof of the property that 
all the roots have the real part equal to $-1/2$ for the Ehrhart polynomials of several polytopes 
including the root polytopes associated with the root systems of type $A$ and $C$.

A root system is a finite set of vectors in $\RR^d$ satisfying certain conditions. 
As is well-known, irreducible root systems are classified as four infinite families of classical root systems of type $A$, $B$, $C$, and $D$, and five exceptional root systems.
For the introduction to root systems, see, e.g., \cite{Hum}. 

A \textit{root polytope} is the convex hull of a root system and we regard it as a lattice polytope with respect to the root lattice. 
In this paper, we mainly discuss the root polytopes of type $A$ and $C$, but the definitions of these root systems, defined later, look different from \cite{Hum}. 
%\cite{ABHPS} and \cite{BHV}. 
%When we use the root systems for types $A$ and $C$,
We apply an appropriate unimodular transformation to make root polytopes as lattice polytopes with respect to the usual lattice $\ZZ^d$. 
Since those are unimodularly equivalent, the information of the lattice polytopes is preserved. In particular, those have the same Ehrhart polynomials.

The classical root system of type $A$ is given by $\{\pm(\bm{e}_i-\bm{e}_j) \; : \; 1 \leq i < j \leq d+1\}$,  
where $\bm{e}_i$ denotes the $i$th coordinate unit vector of $\RR^{d+1}$. 
Note that this set is contained in the hyperplane $H=\{(x_1,\cdots,x_{d+1}) \in \RR^{d+1} : x_1+\cdots+x_{d+1}=0\}$, 
which is isomorphic to $\RR^d$ as vector spaces. Since this set can be transformed into 
\begin{align}\label{typeA}\{\pm(\bm{e}_i+\cdots+\bm{e}_j)\; :\;1\le i \leq j\le d\}\end{align}
via the linear bijection $\bm{e}_i-\bm{e}_{i+1} \mapsto \bm{e}_i \; (i=1,\cdots,d)$ between $H$ and $\RR^d$, 
we use \eqref{typeA} as the root system of type $A$. 
The root polytope of type $A$, denoted by $A_{d}$, is the lattice polytope given as the convex hull of \eqref{typeA}. 
We know that $A_d$ is a reflexive polytope of dimension $d$ and 
the Ehrhart polynomial of the root polytope $A_{d}$ is calculated by Bacher, Harpe, and Venkov \cite[Theorem 3]{BHV} as follows:
\[|kA_{d}\cap\ZZ^d|=\sum_{i=0}^{d}\binom{d}{i}^2\binom{k+d-i}{d}.\]

For the dual polytope $A_{d}^{*}$ of the root polytope $A_{d}$, its Ehrhart polynomial is calculated by Higashitani, Kummer, and Micha\l ek \cite[Lemma 5.3]{HKM} as follows:
\[|kA_{d}^{*}\cap\ZZ^d|=\sum_{i=0}^{d}\binom{d+1}{i}k^i.\]
Furthermore, they proved that both $A_{d}$ and $A_{d}^{*}$ have the property that all the roots of the Ehrhart polynomial have their real parts equal to $-1/2$, 
i.e., those are CL-polytopes. See \cite[Example 3.5 and Corollary 5.4]{HKM}, respectively.

The classical root system of type $C$ is given by $\{\pm (\bm{e}_i \pm \bm{e}_j) : 1 \leq i < j \leq d\} \cup \{\pm 2\bm{e}_i : 1 \leq i \leq d\}$. 
When we take the convex hull of this set, we see that 
$\pm 2 \bm{e}_1,\cdots, \pm 2 \bm{e}_d$ are the vertices and the remaining points are not vertices. 
Hence, $C_d$ is combinatorially a cross polytope, and so $C_d^*$ is a cube. See \cite[Section 4]{ABHPS} for more details. 
Moreover, by applying the linear bijection $\bm{e}_i-\bm{e}_{i+1} \mapsto \bm{e}_i$ ($i=1,\cdots,d-1$) and $2\bm{e}_d \mapsto \bm{e}_d$ 
which sends the root lattice of type $C$ to the usual lattice $\ZZ^d$, we see that 
the set of vectors $\{\pm 2 \bm{e}_i : 1 \leq i \leq d\}$ is transformed into 
\begin{align}\label{typeC}
\{\pm(2\bm{e}_i+\cdots+2\bm{e}_{d-1}+\bm{e}_{d})\; :\;1\le i\le d-1\} \cup \{\pm\bm{e}_d\}. 
\end{align}
Hence, in this paper, we let the root polytope of type $C$, denoted by $C_d$, be the lattice polytope given as the convex hull of \eqref{typeC}. 
Then $C_d$ is a reflexive polytope of dimension $d$. 
The Ehrhart polynomial of $C_d$ is also calculated in \cite[Theorem 3]{BHV} as follows: 
\[|kC_{d}\cap\ZZ^d|=\sum_{i=0}^{d}\binom{2d}{2i}\binom{k+d-i}{d}.\] 
In \cite[Example 3.6]{HKM}, it is proved that $C_d$ is a CL-polytope. 

In this paper, we consider the dual polytope $C_{d}^{*}$ of the root polytope $C_{d}$. We prove the following theorems.

%主定理---------------------------------------------------------------------------
%\section{Main Theorem}
\begin{thm}\label{thm1}
For the dual polytope $C_{d}^{*}$ of the root polytope $C_{d}$, we have 
\[ |k C_{d}^{*}\cap\ZZ^d|=(k+1)^{d}+k^d. \]
\end{thm}

\begin{thm}\label{thm2}
All the roots of the Ehrhart polynomial $E_{C_d^*}(k)$ have the real part equal to $-1/2$. Namely, $C_d^*$ is a CL-polytope. 
Moreover, the series of the Ehrhart polynomials $\{E_{C_d^*}(k)\}_d$ of degree $d$ satisfies the interlacing property. 
\end{thm}

We prove Theorem~\ref{thm1} in Section~\ref{sec:proof1}, and we prove Theorem~\ref{thm2} in Section~\ref{sec:proof2}.

\begin{rem}[{cf. \cite[Remark 3.7]{HKM}}]
For root polytopes of other types, as remarked in \cite[Remark 3.7]{HKM},
the root polytope $B_d$ is not reflexive, hence the roots of the Ehrhart polynomial do not distribute symmetrically with respect to the line $\mathrm{Re}(z)=-1/2$.
The root polytope $D_d$ is reflexive and the roots of the Ehrhart polynomial distribute symmetrically with respect to the line $\mathrm{Re}(z)=-1/2$, however, some of the roots are not on the line.
\end{rem}

\section*{Acknowledgements}
The authors would like to thank Masahiro Hachimori for a lot of his helpful comments on the results. 
The authors would also like to appreciate the numerous comments of the anonymous referees for the previous version of this paper. 
The first named author is partially supported by JSPS Grant-in-Aid for Scientists Research (C) 20K03513.

\medskip

%証明---------------------------------------------------------------------------
\section{The proof of Theorem~\ref{thm1}}\label{sec:proof1}
%In this section, we give the proof of Theorem $1.1$. 
We perform our proof along the following strategy: 
\begin{itemize}
\item[(I)] Reduce the problem to providing a certain bijection between certain sets $\Ac$ and $\Bc$. 
\item[(II)] Define $f : \Ac \rightarrow \Bc$ and $g : \Bc \rightarrow \Ac$. 
\item[(III)] Give the detailed description of $g$ and verify the well-definedness.
\item[(IV)] Prove the well-definedness and the injectivity of $f$. 
\end{itemize}

Note that $f \circ g$ trivially becomes the identity on $\Bc$ by definition, 
and what $g \circ f$ is the identity on $\Ac$ follows from the injectivity of $f$ and the description of $g$.

%Step1----------------------------------------------------------------------------------
\subsection{Step (I): Reduction of the problem} 
For the proof, we use the Ehrhart polynomial of $A_{d}^{*}$. Let us deform the Ehrhart polynomial of $A_{d}^{*}$ as follows: 
\begin{eqnarray*}
|k A_{d}^{*}\cap\ZZ^d|=\sum_{i=0}^{d}\binom{d+1}{i}k^i = \sum_{i=0}^{d+1}\binom{d+1}{i}k^i-k^{d+1} =(k+1)^{d+1}-k^{d+1}.
\end{eqnarray*}
Here, we focus on the boundary of $A_{d}^{*}$. Recall that a lattice polytope $P\subset\RR^d$ is reflexive if and only if 
there is no lattice point in $kP^{\circ}\backslash(k-1)P$ for all positive integers $k$ (\cite{H}), 
i.e. $|kP\backslash(k-1)P\cap\ZZ^d|=|k\partial P\cap\ZZ^d|$ holds, where $\partial X$ denotes the boundary of $X \subset \RR^d$. 
%between $kP$ and $(k-1)P$, 
%for all positive integers $k$ (\cite{H}).
Since $A_d$ is reflexive, $A_d^*$ is also reflexive.
Thus, for $A_{d}^{*}$, the following equation holds: 
\[|k \partial A_{d}^{*}\cap\ZZ^d|=|k A_{d}^{*}\cap\ZZ^d|-|(k-1) A_{d}^{*}\cap\ZZ^d|.\]
Moreover, 
\begin{align*} 
|k \partial A_{d-1}^{*}\cap\ZZ^{d-1}|+2|(k-1) A_{d-1}^{*}\cap\ZZ^{d-1}|
&=|k A_{d-1}^{*}\cap\ZZ^{d-1}|+|(k-1) A_{d-1}^{*}\cap\ZZ^{d-1}| \\
&=(k+1)^d-k^d+k^d-(k-1)^d \\
&=(k+1)^d-(k-1)^d.
\end{align*}

\noindent
{\bf Claim}: For the proof of Theorem~\ref{thm1}, it is enough to show the following equation \eqref{eq:mokuteki} for all positive integers $k$: 
\begin{equation}\label{eq:mokuteki}
|k \partial C_{d}^{*}\cap\ZZ^d|=|k \partial A_{d-1}^{*}\cap\ZZ^{d-1}|+2|(k-1) A_{d-1}^{*}\cap\ZZ^{d-1}|. 
\end{equation}
\begin{proof}
We see that 
\begin{align*}
|k \partial C_{d}^{*}\cap\ZZ^d|&=(k+1)^{d}-(k-1)^{d}, \\
|(k-1) \partial C_{d}^{*}\cap\ZZ^d|&=k^{d}-(k-2)^{d}, \\
|(k-2) \partial C_{d}^{*}\cap\ZZ^d|&=(k-1)^d-(k-3)^d, \\
&\vdots \\
|2 \partial C_{d}^{*}\cap\ZZ^d|&=3^d-1^d, \\
|1 \partial C_{d}^{*}\cap\ZZ^d|&=2^d, \\
|\{0\}\cap\ZZ^d|&=1. 
\end{align*}
Here, we recall that $C_{d}^{*}$ is also reflexive since $C_d$ is reflexive.
Thus, the sum of the left-hand sides of the above equations is \[
|k \partial C_{d}^{*}\cap\ZZ^d|+|(k-1) \partial C_{d}^{*}\cap\ZZ^d|+\cdots+|\{0\}\cap\ZZ^d|=|k C_{d}^{*}\cap\ZZ^d|.\] 
On the other hand, the sum of the right-hand sides is equal to $(k+1)^d+k^d$.
\end{proof}
For the proof of \eqref{eq:mokuteki}, it is enough to show the existence of a bijection 
between $k\partial C_d^*$ and $k\partial A_{d-1}^{*}\sqcup (k-1)A_{d-1}^{*}\sqcup (k-1)A_{d-1}^{*}$. 
Here, the symbol ``$\sqcup$'' stands for a multiset. 
Namely, we provide a bijection from $k \partial C_d^*$ to singly $k\partial A_{d-1}^{*}$ and doubly $(k-1)A_{d-1}^{*}$. 
In other words, we divide $k\partial C_d^*$ into three parts $\Cc,\Cc'$ and $\Cc''$ such that 
there are bijections between the lattice points of $\Cc$ and $k\partial A_{d-1}^*$, 
between those of $\Cc'$ and $(k-1)A_{d-1}^{*}$ and between those of $\Cc''$ and $(k-1)A_{d-1}^{*}$. 
%This bijection means that one essentially divides the polytope $kC_d^*$ into three pieces,
%and that the given maps preserve lattice points. 
Let $$\mathcal{A}=k\partial C_{d}^{*}\;\text{ and }\;\mathcal{B}=k\partial A_{d-1}^{*}\sqcup (k-1)A_{d-1}^{*}\sqcup (k-1)A_{d-1}^{*}.$$ 

\begin{rem}[Comments by referee]
Our algebraic bijection proves the following geometric fact. 
Let $P=k\partial C_d^*$, $Q =kA_{d-1}^*$ and let $\pi : \RR^d \rightarrow \RR^{d-1}$ be the coordinate projection onto the first $d-1$ coordinates, 
which projects the lattice $\ZZ^d$ onto the lattice $\ZZ^{d-1}$. 
Consider the intersection of fibers $\pi^{-1}(x)$ over $x \in Q \cap \ZZ^{d-1}$ with $P$. 
Then $\pi^{-1}(x) \cap P$ contains one lattice point of $P$ if $x \in \partial Q$ 
and two lattice points if $x \in Q \setminus \partial Q$ (one on the lower hull, one on the upper). 
The fibers intersect with $P$ in lattice points as the corresponding coefficients in the constraints of $P$ is $\pm 1$. 
The set of lattice points of the interior of $Q$ is nothing but that of $(k-1)A_{d-1}^*$, so we obtain our bijection. 
\end{rem}

%Step2-----------------------------------------------------------------------------------
\subsection{Step (II): Definitions of $f : \Ac \rightarrow \Bc$ and $g : \Bc \rightarrow \Ac$}
Let $f$ be a map 
\begin{align*}
f : k\partial C_{d}^{*} &\longrightarrow k\partial A_{d-1}^{*}\sqcup (k-1)A_{d-1}^{*}\sqcup (k-1)A_{d-1}^{*}, \\
(\alpha_{1},\cdots,\alpha_{d})&\mapsto(\alpha_{1},\cdots,\alpha_{d-1}). 
\end{align*}
%Here, two copies of $(k-1)A_{d-1}^{*}$ are in the range of $f$. 
%In other view, $f$ maps singly on $k\partial A_{d-1}^{*}$ and doubly on $(k-1)A_{d-1}^{*}$.

Let $g$ be a map 
\begin{align*}
g : k\partial A_{d-1}^{*}\sqcup (k-1)A_{d-1}^{*}\sqcup (k-1)A_{d-1}^{*}&\longrightarrow k\partial C_{d}^{*}, \\
(\alpha_{1},\cdots,\alpha_{d-1})&\mapsto(\alpha_{1},\cdots,\alpha_{d}),
\end{align*}
where $\alpha_{d}$ will be determined suitably below.
%We show that $f$ is bijective.

Each well-definedness of the maps $f$ and $g$ will be proved below. 

\medskip

%%%%
\begin{ex}
We give an example of the bijection when $d=3$ and $k=2$.
All lattice points in $2\partial C_3^*$, $2\partial A_2^*$ and $A_2^*$ are listed in Table~\ref{tab:bijection}. 
From Table~\ref{tab:bijection}, we can check that each lattice point $(\alpha_1,\alpha_2,\alpha_3)$ in $2\partial C_3^* \cap \ZZ^3$ one-to-one corresponds to each lattice point $(\alpha_1,\alpha_2)$ in 
$(2\partial A_2^* \cap \ZZ^2) \sqcup (A_2^* \cap \ZZ^2) \sqcup (A_2^* \cap \ZZ^2)$.

{\rm 
\begin{table}[h]
\centering
\begin{tabular}{|c||c|c|c|} \hline
$2\partial C_3^* \cap \ZZ^3$ & $2\partial A_2^* \cap \ZZ^2$ & $A_2^* \cap \ZZ^2$ & $A_2^* \cap \ZZ^2$\\ \hline\hline
$\pm(2,0,-2)$ & $\pm(2,0)$ &  & \\ \hline
$\pm(2,-1,0)$ & $\pm(2,-1)$ &  & \\ \hline
$\pm(2,-2,2)$ & $\pm(2,-2)$ &  & \\ \hline
$\pm(1,1,-2)$ & $\pm(1,1)$ &  & \\ \hline
$\pm(1,-2,2)$ & $\pm(1,-2)$ &  & \\ \hline
$\pm(0,2,-2)$ & $\pm(0,2)$ &  & \\ \hline
$\pm(0,0,2)$ & &  $(0,0)$ & $(0,0)$ \\ \hline
\begin{tabular}{c}
$\pm(1,0,0)$ \\
$\pm(1,0,-2)$
\end{tabular}
&  & $\pm(1,0)$ & $\pm(1,0)$ \\ \hline
\begin{tabular}{c}
$\pm(1,-1,0)$ \\
$\pm(1,-1,2)$ 
\end{tabular}
 &  & $\pm(1,-1)$ & $\pm(1,-1)$ \\ \hline
\begin{tabular}{c}
$\pm(0,1,0)$ \\ 
$\pm(0,1,-2)$
\end{tabular}
 &  & $\pm(0,1)$ & $\pm(0,1)$ \\ \hline
\end{tabular}
\caption{All lattice points in the case $d=3$ and $k=2$}\label{tab:bijection}
\end{table}
}

\end{ex}

%gとfg=idについて----------------------------------------------------------------------------
%\subsection{The definition of $g$ and the proof of $f\circ g=\mathrm{id}_\mathcal{B}$}\label{sec:fg}
\subsection{Step (III): Details of $g$ and its well-definedness}\label{sec:fg}
%The goal of this section is to give the definition of the map $g$,
%that is, we consider how to construct $\alpha_{d}$ in the map $g$, 
%and check the well-definedness of $g$.  

%\medskip

In order to describe $g$, we divide the domain of $g$ into two parts: $k\partial A_{d-1}^*$ and (two copies of) $(k-1)A_{d-1}^*$. 

%%%%
\subsubsection{}
First, we define the value $\alpha_{d}$ of $g(\alpha_{1},\cdots,\alpha_{d-1})$ for $(\alpha_{1},\cdots,\alpha_{d-1})\in k\partial A_{d-1}^{*}$ 
so that $g(\alpha_{1},\cdots,\alpha_{d-1})$ $\in k\partial C_{d}^{*}$.
By definition of $A_d$,
the polytope $A_d^*$ is the set of $(\alpha_{1},\cdots,\alpha_{d})\in\RR^d$ satisfying the following inequalities: 
\begin{equation}
|\alpha_{i}+\cdots+\alpha_{j}|\leq1\;\;(1\leq i\le j \leq d).\nonumber
\end{equation}
%Note that the above expression is irredundant since \eqref{typeA} precisely coincides with the set of vertices of $A_d$ (see \cite[Section 3]{ABHPS}). 
%(In general, for any reflexive polytope $P$, the facets of $P^*$ one-to-one correspond to the vertices of $P$.) 
Since the point $(\alpha_{1},\cdots,\alpha_{d-1})$ lies on the boundary of $k A_{d-1}^{*}$,
$(\alpha_{1},\cdots,\alpha_{d-1})$ satisfies the following equation:
\begin{eqnarray}\label{eq:1}
|\alpha_{i}+\cdots+\alpha_{j}| = k ~~\mbox{for some}~ i,j ~\mbox{with}~ 1\le i\leq j\le d-1.
\end{eqnarray}

For $(\alpha_{1},\cdots,\alpha_{d-1})\in k \partial A_{d-1}^{*}$,
we define $\alpha_{d}$ of $g(\alpha_{1},\cdots,\alpha_{d-1})$ as follows: \\
%If there exist $i_0$ and $j_0$ such that $|\alpha_{i_0}+\cdots+\alpha_{j_0}|=k$, then let 
%\[\alpha_{d}=-(\alpha_{i_0}+\cdots+\alpha_{d-1})-(\alpha_{j_0+1}+\cdots+\alpha_{d-1}).\]
(i) If there exist $i_0$ and $j_0$ such that 
$\alpha_{i_0}+\cdots+\alpha_{j_0}=k$, then let 
\[\alpha_{d}=-k-2(\alpha_{j_0+1}+\cdots+\alpha_{d-1}).\]
(ii) If there exist $i_0$ and $j_0$ such that
$\alpha_{i_0}+\cdots+\alpha_{j_0}=-k$, then let 
\[\alpha_{d}=k-2(\alpha_{j_0+1}+\cdots+\alpha_{d-1}).\]

\noindent
{\bf Claim}: If there exist two different pairs $(i_0,j_0)$ and $(i_0',j_0')$ with $|\alpha_{i_0}+\cdots+\alpha_{j_0}|=|\alpha_{i_0'}+\cdots+\alpha_{j_0'}|=k$, then the corresponding $\alpha_d$'s are the same. 
\begin{proof}
Let, say, $\alpha_{i_0}+\cdots+\alpha_{j_0}=k$ and $\alpha_{i_0'}+\cdots+\alpha_{j_0'}=k$ with $j_0<j_0'$. Then 
%\begin{align*}
%\overline{\alpha}:=(-k-2(\alpha_{j_0+1}+\cdots+\alpha_{d-1}))-(-k-2(\alpha_{j_0'+1}+\cdots+\alpha_{d-1}))=-2(\alpha_{j_0+1}+\cdots+\alpha_{j_0'}). 
%\end{align*}
\begin{align*}
\alpha_{j_0+1}+\cdots+\alpha_{j_0'}=\begin{cases}
\alpha_{i_0}+\cdots+\alpha_{i_0'-1} &(\text{if }i_0<i_0'), \\
-(\alpha_{i_0'}+\cdots+\alpha_{i_0-1}) &(\text{if }i_0'<i_0), \\
0 &(\text{if }i_0=i_0'). 
\end{cases}
\end{align*}
For our purpose, it suffices to show $\alpha_{j_0+1}+\cdots+\alpha_{j_0'}=0$. 
\begin{itemize}
\item In the case $i_0<i_0'$, we see that 
\begin{align*}
\alpha_{i_0}+\cdots+\alpha_{j_0'}=2k+
\begin{cases}
\alpha_{j_0+1}+\cdots+\alpha_{i_0'-1} &(\text{if }j_0<i_0'-1), \\
-(\alpha_{i_0'}+\cdots+\alpha_{j_0}) &(\text{if }i_0' \leq j_0). 
\end{cases}
\end{align*}
Since $|\alpha_i+\cdots+\alpha_j| \leq k$ holds for any $1 \leq i,j \leq d-1$, we obtain that 
\begin{align*}
\begin{cases}
\alpha_{j_0+1}+\cdots+\alpha_{i_0'-1}=-k &(\text{if }j_0<i_0'-1), \\
\alpha_{i_0'}+\cdots+\alpha_{j_0}=k &(\text{if }i_0' \leq j_0). 
\end{cases}
\end{align*}
Hence, we obtain that \begin{align*}
\alpha_{j_0+1}+\cdots+\alpha_{j_0'}=\begin{cases}
(\alpha_{j_0+1}+\cdots+\alpha_{i_0'-1})+(\alpha_{i_0'}+\cdots+\alpha_{j_0'}) &(\text{if }j_0<i_0'-1) \\
(\alpha_{i_0'}+\cdots+\alpha_{j_0'})-(\alpha_{i_0'}+\cdots+\alpha_{j_0}) &(\text{if }i_0' \leq j_0) 
\end{cases}
=0. \end{align*}
\item In the case $i_0'<i_0$, by \eqref{eq:7} below, we know that $\alpha_{j_0+1}+\cdots+\alpha_{j_0'} \leq 0$, 
while we also know that $\alpha_{j_0+1}+\cdots+\alpha_{j_0'}=-(\alpha_{i_0'}+\cdots+\alpha_{i_0-1}) \geq 0$, as required.
\end{itemize}
The other cases (e.g., $\alpha_{i_0}+\cdots+\alpha_{j_0}=k$ and $\alpha_{i_0'}+\cdots+\alpha_{j_0'}=-k$) can be also discussed in the similar way. 
\end{proof}

%\begin{align*}
%|\alpha_{i}|&\leq 1\;\; (1\le i\le d),\\
%|\alpha_{i}+\cdots+\alpha_{j}|&\leq 1\;\; (1\le i<j\le d-1), \mbox{ and} \\
%|2(\alpha_{i}+\cdots+\alpha_{d-1})+\alpha_{d})|&\leq 1\;\; (1\le i\le d-1).
%\end{align*}

Now, we verify the well-definedness of $g$, i.e., $(\alpha_{1},\cdots,\alpha_{d})$ lies in $k\partial C_{d}^{*}$. 
By definition of $C_d$, we see that $(\alpha_1,\cdots,\alpha_d) \in kC_d^*$ if and only if it satisfies the following inequalities: 
\begin{align}
|\alpha_{d}| &\leq k, \text{ and} \label{eq:5}\\
|2(\alpha_i+\cdots+\alpha_{d-1})+\alpha_{d}| &\le k \;\;(1\leq i\leq d-1). \label{eq:6}
\end{align}
%We assume the case (i). 
We have to check that these inequalities hold in the cases (i) and (ii), but we only prove for the case (i) since the case (ii) can be proved similarly.
Note that since $(\alpha_1,\cdots,\alpha_{d-1}) \in kA_{d-1}^*$, we know that the following inequalities hold: 
\begin{equation}\label{eq:3}
|\alpha_{i}+\cdots+\alpha_{j}| \le k ~~(1\le i\le j\le d-1).
\end{equation}

We can see that the points on $k\partial C_{d}^{*}$ satisfy the following three inequalities: 
%For the case (i), 
%we use that the points on $k\partial C_{d}^{*}$ satisfy the following three inequalities.
%we use the following three inequalities: 
\begin{align}\label{eq:7} 
\begin{cases}
\alpha_{l}+\cdots+\alpha_{i_0-1}\le0 & \;\;(l\le i_0-1),\\
\alpha_{i_0}+\cdots+\alpha_{l}\ge0 & \;\;(i_0\le l\le j_0),\\
\alpha_{j_0+1}+\cdots+\alpha_{l}\le0 & \;\;(j_0+1\le l).
\end{cases}
\end{align}
%We use them for verifying \eqref{eq:5} and \eqref{eq:6}.
In fact, if there exists $l$ with $\alpha_{l}+\cdots+\alpha_{i_0-1}>0$, 
then $\alpha_l+\cdots+\alpha_{i_0-1}+\alpha_{i_0}+\cdots+\alpha_{j_0}>0+k$, a contradiction to $(\alpha_1,\cdots,\alpha_{d-1}) \in k\partial A_{d-1}^{*}$. 
We can verify the other equalities in the same way. 

Here, from \eqref{eq:3} and \eqref{eq:7}, we have $-k\leq \alpha_{j_0+1}+\cdots+\alpha_{d-1}\leq0$, which implies that 
$-k\leq -k-2(\alpha_{j_0+1}+\cdots+\alpha_{d-1})\leq k$. Therefore, we have \[-k\leq\alpha_{d}\leq k,\] 
and \eqref{eq:5} is verified.
%from \eqref{eq:6} and the first statement of the inequalities $(9)$, 
%\[\alpha_{j+1}+\cdots+\alpha_{d-1}=-k.\]
%Therefore $-k\leq \alpha_{d}=-k-2(\alpha_{j+1}+\cdots+\alpha_{d-1})\leq k$, Hence (7) is verified. 
On the other hand, for $l\le i_0-1$, we have 
%Consider $\alpha_{d}=-k-2(\alpha_{j+1}+\cdots+\alpha_{d-1})$. 
\begin{align*}
2(&\alpha_{l}+\cdots+\alpha_{d-1})+\alpha_{d} \\
&=2(\alpha_{l}+\cdots+\alpha_{i_0-1})+2k+2(\alpha_{j_0+1}+\cdots+\alpha_{d-1})-k-2(\alpha_{j_0+1}+\cdots+\alpha_{d-1}) \\
&=k+2(\alpha_{l}+\cdots+\alpha_{i_0-1}), 
\end{align*}
and from \eqref{eq:3} and \eqref{eq:7}, we have $-k\leq\alpha_{l}+\cdots+\alpha_{i_0-1}\leq0$, which implies that 
$-2k\leq2(\alpha_{l}+\cdots+\alpha_{i_0-1})\leq0$. 
Therefore, we have 
\[-k\leq k+2(\alpha_{j_0+1}+\cdots+\alpha_{d-1})\leq k,\]
and \eqref{eq:6} is verified. 
For the cases that $i_0\leq l\leq j_0$ and $j_0+1\leq l$, \eqref{eq:6} is similarly verified. 
Hence, we conclude $(\alpha_{1},\cdots,\alpha_{d})\in k\partial C_{d}^{*}$.

The discussion for the case (ii) is similar.
%Also, for the case (ii) we apply the same way. 
In this case, we use that the points on $k\partial C_{d}^{*}$ satisfy
the following three inequalities:
%we suppose $\alpha_{i}+\cdots+\alpha_{j}=-k$, and put $\alpha_{d}=k-2(\alpha_{j+1}+\cdots+\alpha_{d-1})$.
%And for $(\alpha_{1},\cdots,\alpha_{d-1})$ the follow inequalities hold.
\[\begin{cases}
\alpha_{l}+\cdots+\alpha_{i_0-1}\ge0 & (l\le i_0-1),\\
\alpha_{i_0}+\cdots+\alpha_{l}\le0 & (i_0\le l\le j_0),\\
\alpha_{j_0+1}+\cdots+\alpha_{l}\ge0 & (j_0+1\leq l).
\end{cases}
\]
We can verify \eqref{eq:5} and \eqref{eq:6} from these inequalities, 
and we conclude $(\alpha_{1},\cdots,\alpha_{d})\in k\partial C_{d}^{*}$. 

\medskip

%%%%
\subsubsection{}
Next, we define $\alpha_{d}$ of $g(\alpha_{1},\cdots,\alpha_{d-1})$ for $(\alpha_{1},\cdots,\alpha_{d-1}) \in (k-1)A_{d-1}^{*}$ in two ways 
so that $g(\alpha_{1},\cdots,\alpha_{d-1}) \in k\partial C_{d}^{*}$.
%Next, we verify $g(\alpha_{1},\cdots,\alpha_{d-1})\in k \partial C_{d}^{*}$ for $(\alpha_{1},\cdots,\alpha_{d-1})\in(k-1)A_{d-1}^{*}\sqcup(k-1)A_{d-1}^{*}$.
%the part $(k-1)A_{d-1}^{*}\sqcup(k-1)A_{d-1}^{*}$.
From the definition of $A_{d}^{*}$, all $(\alpha_{1},\cdots,\alpha_{d-1})\in(k-1)A_{d-1}^{*}$ satisfy the following inequalities:
%\[|\alpha_{i}|\le k ~~~~~(1\le i\le d-1)\]
\begin{eqnarray}\label{eq:8}
|\alpha_{i}+\cdots+\alpha_{j}| \le k-1 ~~(1\le i\le j\le d-1).
\end{eqnarray}
Here, we set 
\begin{align*}
&p:=\max\{\alpha_{i}+\cdots+\alpha_{d-1} : 1\le i\le d-1\}\text{ and }\\
&q:=\min\{\alpha_{i}+\cdots+\alpha_{d-1} : 1\le i\le d-1\}.
\end{align*}
Note that $-k+1\le p\le k-1$, $-k+1\le q\le k-1$, and $p-q\le k-1$.
We set two types of $\alpha_{d}$ as follows:
\begin{align}\label{eq:cases}
\alpha_{d}^{(1)}=
\begin{cases}
k,& \text{if }p\le0,\\
k-2p,& \text{if }0<p,\\
\end{cases}
~~~~~~~~{\rm and}~~~~~~~\alpha_{d}^{(2)}=
\begin{cases}
-k-2q, & \text{if }q<0,\\
-k, & \text{if }0\le q. \end{cases}\end{align}

\noindent
{\bf Claim}: $\alpha_{d}^{(1)} \neq \alpha_{d}^{(2)}$. 
\begin{proof}
For example, we consider the situations of $0<p$ and $q<0$, and assume that $\alpha_{d}^{(1)}=\alpha_{d}^{(2)}$. 
Then we have $k-2p=-k-2q$, so we have $k=p-q\le k-1$, a contradiction. 
For the other cases, we can check in the same way. Therefore, $\alpha_{d}^{(1)}$ and $\alpha_{d}^{(2)}$ are different. 
%and there are only two cases $\alpha_{d}^{(1)}$ and $\alpha_{d}^{(2)}$.
\end{proof}
Hence, each $(\alpha_{1},\cdots,\alpha_{d-1})\in(k-1)A_{d-1}^{*}$ has two different $\alpha_{d}$'s ($\alpha_{d}^{(1)}$ and $\alpha_{d}^{(2)}$).

Now, we verify that each $(\alpha_{1},\cdots,\alpha_{d})$ satisfies the definitions of $k\partial C_{d}^{*}$.
%What we want to show is that one of the three inequalities \eqref{eq:2},\eqref{eq:3} and \eqref{eq:4} holds with equality. 
%since $(\alpha_{1},\cdots,\alpha_{d})$ lies on the boundary of $kC_{d}^{*}$.
%From \eqref{eq:8} the inequality \eqref{eq:3} holds with strict inequality and the inequality \eqref{eq:2} holds with strict inequality for $1\leq i\leq d-1$. 
We need to show that \eqref{eq:5} and \eqref{eq:6} are true and one equality in \eqref{eq:5} and \eqref{eq:6} holds.
Namely, we check that each $(\alpha_{1},\cdots,\alpha_{d})$ satisfies either  
\begin{align}
|\alpha_{d}|&= k \text{ or } \label{eq:9} \\
|2(\alpha_i+\cdots+\alpha_{d-1})+\alpha_{d}|&= k. \label{eq:10}
\end{align}
\begin{itemize}
\item When $\alpha_{d}^{(1)}=k~(p\leq0)$, \eqref{eq:9} clearly holds.
On the other hand, we have 
%So for each $(\alpha_{1},\cdots,\alpha_{d})$ we check the following statements are correct:\\
%\[
%\begin{cases}
%\cdot~\mbox{If the case (i) is satisfied, then}~~|2(\alpha_{i}+\cdots+\alpha_{d-1})+\alpha_{d}|\le k.\\
%\cdot~\mbox{If the case (ii) is satisfied, then}~~|\alpha_{d}|\le k.
%\end{cases}
%\]
%If the case (i) is satisfied, then
%\begin{equation}
%|2(\alpha_{i}+\cdots+\alpha_{d-1})+\alpha_{d}|\le k,
%\end{equation}
%and if the case (ii) is satisfied, then
%\begin{equation}
%|\alpha_{d}|\le k.
%\end{equation}
%We consider case (i). Since if $\alpha_{d}=k$, then $p\leq 0$, 
\[-k<2(-k+1)+k \leq 2(\alpha_{i}+\cdots+\alpha_{d-1})+\alpha_{d}\le2p+k\leq k.\]
Hence, \eqref{eq:6} is verified. 
\item When $\alpha_{d}^{(1)}=k-2p~(0<p)$, \eqref{eq:10} clearly holds by choosing $\max\{\alpha_{i}+\cdots+\alpha_{d-1}\}$. 
On the other hand, we have \[-k<k-2(k-1) \leq k-2p=\alpha_d <k.\]
Hence, \eqref{eq:5} holds.
\item When $\alpha_{d}^{(2)}=-k-2q~(q<0)$, \eqref{eq:10} clearly holds by choosing $\min\{\alpha_{i}+\cdots+\alpha_{d-1}\}$. 
On the other hand, we have \[-k<-k-2q=\alpha_d \leq -k-2(-k+1)<k.\] Hence, \eqref{eq:5} holds.
\item When $\alpha_{d}^{(2)}=-k~(0\leq q)$, \eqref{eq:9} clearly holds. On the other hand, we have
\[-k\leq 2q - k \leq 2(\alpha_{i}+\cdots+\alpha_{d-1})+\alpha_{d}\le 2(k-1)-k<k.\] 
Hence, \eqref{eq:6} is also verified.
\end{itemize}
Therefore, we conclude $(\alpha_{1},\cdots,\alpha_{d})\in k\partial C_{d}^{*}$.
%%%%%

%Finally, by the construction of $g$, we see that the map $g$ is injective.

%gf=id------------------------------------------------------------------------
\subsection{Step (IV): The well-definedness and the injectivity of $f$}
The goal of this section is to check the well-definedness and the injectivity of $f$. 
Recall that the map $f$ just forgets the last entry. 

Consider $(\alpha_{1},\cdots,\alpha_{d})\in k\partial C_{d}^{*}$.
Since it lies on the boundary, one of the inequalities \eqref{eq:5} and \eqref{eq:6} holds with equality. 
For the well-definedness of $f$, we have to check that $(\alpha_1,\cdots,\alpha_{d-1})$ satisfies \eqref{eq:3}, 
but this can be easily obtained from \eqref{eq:5} and \eqref{eq:6} as follows: for any $1 \leq i \leq j \leq d-1$, we see that 
\begin{align*}
|\alpha_i+\cdots+\alpha_j| &= \frac{1}{2}|2(\alpha_i+\cdots+\alpha_{d-1})+\alpha_d-(2(\alpha_{j+1}+\cdots+\alpha_{d-1})+\alpha_d)| \\
&\leq \frac{1}{2}(|2(\alpha_i+\cdots+\alpha_{d-1})+\alpha_d|+|2(\alpha_{j+1}+\cdots+\alpha_{d-1})+\alpha_d|) \\
&\leq k, 
\end{align*}
where we regard $2(\alpha_{j+1}+\cdots+\alpha_{d-1})+\alpha_d$ as $\alpha_d$ when $j=d-1$.  

In what follows, we prove the injectivity of $f$. 

%%%%%
\subsubsection{}
First, we consider the case of $(\alpha_{1},\cdots,\alpha_{d-1})\in k\partial A_{d-1}^{*}$. From the definition of $k\partial A_{d-1}^{*}$, there exist $i_0$ and $j_0$ $(1\leq i_0\leq j_0\leq d-1)$ such that
\begin{eqnarray}%\label{eq:11}
|\alpha_{i_0}+\cdots+\alpha_{j_0}|= k.
\end{eqnarray}
%Note that from the definition of $k\partial C_{d}^{*}$, the inequalities \eqref{eq:5} and \eqref{eq:6} hold at the same time.
\begin{itemize}
\item Assume that $i_0$ and $j_0$ satisfy $\alpha_{i_0}+\cdots+\alpha_{j_0}=k$. %in \eqref{eq:11}.
From \eqref{eq:6}, we see the following: 
\begin{align*}
2(\alpha_{i_0}+\cdots+\alpha_{d-1})+\alpha_{d}&\leq k, \\
\alpha_{d}&\leq  k-2(\alpha_{i_0}+\cdots+\alpha_{j_0})-2(\alpha_{j_0+1}+\cdots+\alpha_{d-1}), \\
\alpha_{d}&\leq -k-2(\alpha_{j_0+1}+\cdots+\alpha_{d-1}).\nonumber
\end{align*}
On the other hand, we also have 
\begin{align*}
-k&\leq 2(\alpha_{j_0+1}+\cdots+\alpha_{d-1})+\alpha_{d},\\
-k-2(\alpha_{j_0+1}+\cdots+\alpha_{d-1})&\leq \alpha_{d}. 
\end{align*}
Thus, $\alpha_{d}=-k-2(\alpha_{j_0+1}+\cdots+\alpha_{d-1})$.
%Therefore, 
%\[\alpha_{d}=-k-2(\alpha_{j+1}+\cdots+\alpha_{d-1}).\]
\item In the case $\alpha_{i_0}+\cdots+\alpha_{j_0}=-k$, %in \eqref{eq:11}, 
from \eqref{eq:6}, we have 
\begin{align*}
-k &\leq 2(\alpha_{i_0}+\cdots+\alpha_{d-1})+\alpha_{d}, \\
-k-2(\alpha_{i_0}+\cdots+\alpha_{j_0})-2(\alpha_{j_0+1}+\cdots+\alpha_{d-1}) &\leq \alpha_{d}, \\
k-2(\alpha_{j_0+1}+\cdots+\alpha_{d-1})&\leq\alpha_{d}. 
\end{align*}
On the other hand, we also have 
\begin{align*}
2(\alpha_{j_0+1}+\cdots+\alpha_{d-1})+\alpha_{d}&\leq k, \\
\alpha_{d}&\leq k-2(\alpha_{j_0+1}+\cdots+\alpha_{d-1}). 
\end{align*}
Thus, $\alpha_{d}=k-2(\alpha_{j_0+1}+\cdots+\alpha_{d-1})$.
\end{itemize}

These discussions imply that $f$ is injective if $f(\alpha_1,\cdots,\alpha_d) \in k\partial A_{d-1}^*$. 

\medskip

%%%%
\subsubsection{}
Next, we consider the case of $(\alpha_{1},\cdots,\alpha_{d-1})\in (k-1)A_{d-1}^{*}$.
%From the definition of $(k-1)A_{d-1}^{*}$, there exist $i_0$ and $j_0$ such that $|\alpha_{i_0}+\cdots+\alpha_{j_0}|\leq k-1$. 
%Thus, \eqref{eq:2} and \eqref{eq:3} are satisfied but not satisfied with equality for $1\leq i\leq d-1$. 
Since $(\alpha_1,\cdots,\alpha_d) \in k \partial C_d^*$, either $|\alpha_{d}|=k$ or $|2(\alpha_{l}+\cdots+\alpha_{d-1})+\alpha_{d}|=k$ holds for some $1 \leq l \leq d-1$. 
In what follows, we use the same notation on $p$ and $q$ as in Subsection~\ref{sec:fg}. 
\begin{itemize}
\item Let $\alpha_{d}=k$. By \eqref{eq:6}, we have $2p+k \leq k$, i.e., $p \leq 0$. 
Hence, we have $$\max\{\alpha_{i}+\cdots+\alpha_{d-1}:1 \leq i \leq d-1\}\leq0.$$ \item Let $\alpha_{d}=-k$. We can show $q \geq 0$ in the same way. 
\item Let $2(\alpha_{l}+\cdots+\alpha_{d-1})+\alpha_{d}=k$. 
If there exists $l'$ with $\alpha_{l}+\cdots+\alpha_{d-1}<\alpha_{l'}+\cdots+\alpha_{d-1}$, 
since $2(\alpha_{l'}+\cdots+\alpha_{d-1})+\alpha_{d}\leq k$ from \eqref{eq:6}, 
we have $2(\alpha_{l'}+\cdots+\alpha_{d-1})+\alpha_{d} \leq 2(\alpha_{l}+\cdots+\alpha_{d-1})+\alpha_{d}$, a contradiction. 
Thus, \[\alpha_{l}+\cdots+\alpha_{d-1}=\max\{\alpha_{i}+\cdots+\alpha_{d-1}:1 \leq i \leq d-1\}=p.\]
Moreover, we have $k-2p \leq k$, i.e., $p \geq 0$. 
In particular, since $\alpha_{d}=k$ when $p=0$, we only need to assume $p>0$. 
%\item Let $2(\alpha_{l}+\cdots+\alpha_{d-1})+\alpha_{d}=-k$. Then we can show $\alpha_{d}=-k-2q$ and $q<0$ in the same way. 
\item In the case where $2(\alpha_{l}+\cdots+\alpha_{d-1})+\alpha_{d}=-k$, we can similarly show that $q<0$. 
\end{itemize}
By these discussions, we can reconstruct $\alpha_d$ in two ways and those are different (see \eqref{eq:cases} and the claim below). 
This means that $f$ is injective if $f(\alpha_1,\cdots,\alpha_d) \in (k-1)A_{d-1}^*$ since $(k-1)A_{d-1}^*$ is prepared as the codomain of $f$. 

\medskip

Therefore, we conclude that the map $f$ is injective.

%This completes the proof of Theorem~\ref{thm1}

\medskip

%CL性interlace性----------------------------------------------------------------------------
\section{The proof of Theorem~\ref{thm2}}\label{sec:proof2}
Let $\alpha\in\CC$ be a root of the Ehrhart polynomial $|kC_{d}^{*}\cap\ZZ^d|=(k+1)^d+k^d$. Then we have
\begin{align*}
\alpha^d=-(\alpha+1)^d \;\Longleftrightarrow\; \left(\frac{\alpha}{\alpha+1}\right)^d=-1 
\;\Longleftrightarrow\;\frac{\alpha}{\alpha+1}=e^{i\frac{2k-1}{d}\pi}~~~(k=1,\cdots,d~). 
\end{align*}
By setting $r=e^{i\frac{2k-1}{d}\pi}$, we have
\[\alpha=\frac{r}{1-r}.\]
On the other hand, we have $r=\cos\theta_{(d,k)}+i\sin\theta_{(d,k)}$, where we let $\theta_{(d,k)}=\frac{2k-1}{d}\pi$. Hence, 
\begin{align*}
\alpha%=\frac{r}{1-r}
=\frac{\cos\theta_{(d,k)}+i\sin\theta_{(d,k)}}{1-\cos\theta_{(d,k)}-i\sin\theta_{(d,k)}} 
=\frac{\cos\theta_{(d,k)}-1+i\sin\theta_{(d,k)}}{2(1-\cos\theta_{(d,k)})} =-\frac{1}{2}+\frac{\sin\theta_{(d,k)}}{2(1-\cos\theta_{(d,k)})}i. 
\end{align*}
This completes the proof of the first argument from Theorem~\ref{thm2}.

Moreover, since 
\[\frac{d}{d\theta}\frac{\sin\theta}{2(1-\cos\theta)}=\frac{\cos\theta-1}{2(1-\cos\theta)^2}<0 \;\; (0<\theta<2\pi),\] 
we see that the imaginary part of the roots is a monotonic decreasing function. 
%So to prove the imaginary part is interlacing, for $\theta_{(d,k)}$ and $\theta_{(d+1,k)}$ 
%it is enough to show the magnitude relation of the angle as follows
Hence, the interlacing property directly follows from the following relations: 
\[\frac{1}{d+1}\pi<\frac{1}{d}\pi<\frac{3}{d+1}\pi<\cdots<\frac{2d-1}{d+1}\pi<\frac{2d-1}{d}\pi<\frac{2d+1}{d+1}\pi.\]
%In fact, we can check easily
%\[\frac{2k-1}{d+1}\pi\le\frac{2k-1}{d}\pi\]
%and
%\[\frac{2k-1}{d}\pi\le\frac{2k+1}{d+1}\pi.\]
%Therefore, $\theta_{(d,k)}$ is interlacing polynomial.
%Therefore, the imaginary part is a interlacing polynomial.

This completes the proof of the second argument from Theorem~\ref{thm2}.

\end{document}